\documentclass[12pt]{article}

\newcommand\eq[1] {(\ref{#1})}

\newcommand{\dund}[1]{\underline{\underline{#1}}}

\newcommand{\bfm}[1]{\mbox{\boldmath ${#1}$}}
\newcommand{\nonum}{\nonumber \\}
\newcommand{\beqa}{\begin{eqnarray}}
\newcommand{\eeqa}[1]{\label{#1}\end{eqnarray}}
\newcommand{\beq}{\begin{equation}}
\newcommand{\eeq}[1]{\label{#1}\end{equation}}
\newcommand{\Grad}{\nabla}
\newcommand{\Div}{\nabla \cdot}

\newcommand{\Real}{\mathop{\rm Re}\nolimits}
\newcommand{\Imag}{\mathop{\rm Im}\nolimits}
\newcommand{\Tr}{\mathop{\rm Tr}\nolimits}

\newcommand{\lang}{\langle}
\newcommand{\rang}{\rangle}
\newcommand{\Md}{\partial}

\newcommand{\Ga}{\alpha}
\newcommand{\Gb}{\beta}
\newcommand{\Gd}{\delta}
\newcommand{\Ge}{\epsilon}

\newcommand{\Gg}{\gamma}

\newcommand{\GD}{\Delta}

\newcommand{\GO}{\Omega}

\newcommand{\BGD}{\bfm\Delta}

\newcommand{\BGG}{\bfm\Gamma}

\newcommand{\CC}{{\cal C}}

\newcommand{\CE}{{\cal E}}

\newcommand{\CJ}{{\cal J}}
\newcommand{\CK}{{\cal K}}

\newcommand{\CM}{{\cal M}}

\newcommand{\CS}{{\cal S}}
\newcommand{\CT}{{\cal T}}
\newcommand{\CU}{{\cal U}}

\def\Ba{{\bf a}}

\def\Bk{{\bf k}}

\def\Bn{{\bf n}}

\def\Bq{{\bf q}}

\def\Bt{{\bf t}}

\def\Bv{{\bf v}}
\def\Bw{{\bf w}}
\def\Bx{{\bf x}}

\def\BA{{\bf A}}

\def\BD{{\bf D}}
\def\BE{{\bf E}}

\def\BG{{\bf G}}
\def\BH{{\bf H}}
\def\BI{{\bf I}}
\def\BJ{{\bf J}}

\def\BL{{\bf L}}
\def\BM{{\bf M}}

\def\BP{{\bf P}}

\def\BR{{\bf R}}
\def\BS{{\bf S}}
\def\BT{{\bf T}}
\def\BU{{\bf U}}
\def\BV{{\bf V}}
\def\BW{{\bf W}}

\def \RR {{\mathbb R}}
\def \CC {{\mathbb C}}

\def \ba {\begin{array}}
\def \ea {\end{array}}
\newtheorem {Thm} {Theorem} [section]

\newtheorem {Adef} [Thm] {Definition}

\newtheorem {Arem} [Thm] {Remark}

\newtheorem {Aexa} [Thm] {Example}

\newtheorem {Anot} [Thm] {Notation}

\def \refe #1.{(\ref{#1})}
\def \reff #1.{figure~\ref{#1}}
\def \refs #1.{section~\ref{#1}}
\def \refss #1.{subsection~\ref{#1}}
\def \refD #1.{Definition~\ref{#1}}
\def \refT #1.{Theorem~\ref{#1}}
\def \refL #1.{Lemma~\ref{#1}}
\def \refC #1.{Corollary~\ref{#1}}
\def \refP #1.{Proposition~\ref{#1}}
\def \refR #1.{Remark~\ref{#1}}
\def \refE #1.{Example~\ref{#1}}
\def \refN #1.{Notation~\ref{#1}}

\usepackage{graphics}
\usepackage{amssymb}
\begin{document}
\vspace{-1in}
\title{Sharp inequalities which generalize the divergence theorem--an extension of the notion of quasiconvexity}

\author{Graeme W. Milton\\
\small{Department of Mathematics, University of Utah, Salt Lake City UT 84112, USA} \\
\small{(milton@math.utah.edu)}}
\date{}
\maketitle
\begin{abstract}

Subject to suitable boundary conditions being imposed, sharp inequalities are obtained on integrals over a region $\GO$ of certain special quadratic functions $f(\BE)$ where $\BE(\Bx)$ derives from a potential $\BU(\Bx)$. With $\BE=\Grad\BU$ it is known that such sharp inequalities can be obtained  when $f(\BE)$ is a quasiconvex function and when $\BU$ satisfies affine boundary conditions (i.e., for some matrix $\BD$, $\BU=\BD\Bx$ on $\Md\GO$). Here we allow
for other boundary conditions and for fields $\BE$ that involve derivatives of a variety orders of $\BU$. We
define a notion of convexity that generalizes quasiconvexity. $Q^*$-convex quadratic functions are introduced, characterized
and an algorithm is given for generating sharply $Q^*$-convex functions. We emphasize that this also solves the outstanding problem of
finding an algorithm for generating extremal quasiconvex quadratic functions. We also treat integrals over $\GO$ of special
quadratic functions $g(\BJ)$ where $\BJ(\Bx)$ satisfies a differential constraint involving derivatives with, possibly, a variety of orders. The results generalize an example of Kang and the author in three spatial dimensions where $\BJ(\Bx)$ is a $3\times 3$ matrix valued field satisfying $\Div\BJ=0$. 

\end{abstract}
\vskip2mm

\noindent Keywords: Quasiconvexity, $Q^*$-convexity, Null-Lagrangians 
\section{Introduction}
\setcounter{equation}{0}
The divergence theorem
\beq \int_{\GO}\Div \BU(\Bx)\,d\Bx=\int_{\partial\GO}\BU\cdot\Bn\,dS,
\eeq{0.1}
of Lagrange, Gauss , Ostrogradsky, and Green, is one of the most important theorems in mathematics, particularly in applied mathematics. 
(Here $\Bn$ is the outward normal to the surface $\Md\GO$ of the region $\GO$). In one dimension it reduces to the fundamental
theorem of calculus and by setting $\BU=w\BV$ one obtains the multidimensional version of integration by parts:
\beq \int_{\GO}\BV(\Bx)\cdot\Grad w\,d\Bx =\int_{\partial\GO}w\BV\cdot\Bn\,dS-\int_{\GO}w\Div \BV(\Bx)\,d\Bx,
\eeq{0.2}
which with $\BV=\Grad v$ yields Green's first identity,
\beq \int_{\GO}\Grad v \cdot\Grad w\,d\Bx =\int_{\partial\GO}w\Grad v\cdot\Bn\,dS-\int_{\GO}w\GD v(\Bx)\,d\Bx.
\eeq{0.3}
Note that the left hand side of \eq{0.1} can be reexpressed as
\beq \int_{\GO}\Tr\BE\,d\Bx, \eeq{0.4}
where $\BE=\Grad \BU$ is the gradient of the potential $\BU$. 

This immediately leads to the question: for what functions $f$ can one express 
\beq \int_\GO f(\BE(\Bx))\,d{\Bx} \eeq{0.5}
in terms of boundary data, when $\BE$ derives from a potential? Such functions are known as null-Lagrangians. One example in two-dimensions 
of a null-Lagrangian is the determinant of $\BE$ with $\BE=\Grad\BU$ and $\BU$ being a two-component vector: one has
\beq f(\BE)=\det\BE=\frac{\Md U_1}{\Md x_1}\frac{\Md U_2}{\Md x_2}-\frac{\Md U_1}{\Md x_2}\frac{\Md U_2}{\Md x_1}=
\frac{\Md}{\Md x_1}\left[U_1\frac{\Md U_2}{\Md x_2}\right]+\frac{\Md}{\Md x_2}\left[-U_1\frac{\Md U_2}{\Md x_1}\right],
\eeq{0.6}
and we directly see that the quantity on the right is a divergence. With
\beq \BE=\pmatrix{\BE_1 \cr \BE_2}=\pmatrix{\Grad \BU \cr \BU},
\eeq{0.7}
another null-Lagrangian is 
\beqa f(\BE)=(\BE_2\cdot\BE_2)\Tr(\BE_1)+2\BE_2\cdot[(\BE_1\BE_2)] & = & (\BU\cdot\BU)\Div\BU+2\BU\cdot[(\Grad\BU)\BU] \nonum
                                                                                                       & = & \Div[\BU(\BU\cdot\BU)],
\eeqa{0.8}
which again is a divergence. 
It has been shown by Ball, Currie and Olver \cite{Ball:1981:NLW} that any $C^1$ null-Lagrangian can be 
expressed as a divergence, so that the evaluation of \eq{0.5} reduces to an application of the divergence theorem. Necessary and essentially sufficient algebraic
conditions to determine whether a (quadratic or nonquadratic) function is a null Lagrangian have been given by 
Murat \cite{Murat:1978:CPC, Murat:1981:CPC, Murat:1987:SCC} (see also \cite{Pedregal:1989:WCW}). In the important case where $\BE=\Grad\BU$, $f(\BE)$ is a null Lagrangian if and only if it is a linear combination of the
subdeterminants (minors) of any order of $\BE$. For references see the paper of Ball, Currie, and Olver, who also show that when $\BE(\Bx)=\Grad^k\BU(\Bx)$ there are no new null Lagrangians $f(\BE)$ beyond those obtained by applying
the result for $k=1$ to the $d\ell^{k-1}$-component potential $\Grad^{k-1}\BU(\Bx)$.

Instead of seeking equalities we search for sharp inequalities and try to find functions $f$ and boundary data for which one can obtain
sharp bounds on integrals of the form \eq{0.5}. Thus, whereas the divergence theorem is an equality relating to the integral of a linear function of $\BE$ to boundary fields, we derive sharp inequalities relating the integral of certain quadratic functions of $\BE$ to boundary fields, for certain boundary fields. This is the sense in which our sharp inequalities generalize the divergence theorem.   If $\BE=\Grad\BU$ and affine boundary conditions $ \BU=\BD\Bx$ on $\Md\GO$ are imposed, then one
possible value of $\BE$ in the interior of $\GO$ is of course $\BE=\BD$ (we follow the usual convention that $\Grad U_1, \Grad U_2, \ldots{},\Grad U_{\ell}$ are the 
rows of $\Grad\BU$, so that $\BE$ is a $\ell\times d$ matrix), and functions $f$ that satisfy the inequality
\beq  \int_\GO f(\BE(\Bx))\,d{\Bx}\geq  \int_\GO f(\BD)\,d{\Bx}=|\GO|f(\BD)
\eeq{0.9}
are called quasiconvex and for such functions the inequality is obviously sharp. (Here $|\GO|$ denotes the volume of $\GO$). For a good introduction to quasiconvexity the reader is referred to the book of Dacorogna \cite{Dacorogna:2007:DMC} and references therein. As discussed there, examples of quasiconvex functions include convex
and polyconvex functions. So far only the quadratic quasiconvex functions have been completely characterized: see 
Tartar and Murat \cite{Tartar:1979:CCA, Murat:1985:CVH, Tartar:1985:EFC}. In particular, the quadratic function
$f(\Grad\BU)$ is quasiconvex if and only if $f(\BH)$ is non-negative for all  rank one $\ell\times d$ matrices $\BH$.
One of contributions of this paper is to show that for quadratic quasiconvex functions that are sharply
quasiconvex, sharp bounds on the integral \eq{0.5} can be obtained for a wide variety of boundary conditions, and not just affine ones. We 
emphasize (see p.26 of \cite{Morrey:1952:QSM} and \cite{Ball:1985:BPI}) that it follows
 from an example of Terpstra\cite{Terpstra:1938:DBF}, and was shown more simply  by Serre\cite{Serre:1981:CLH, Serre:1983:ROE},
that there are quadratic quasiconvex functions that are not the sum of a convex quadratic function and a null-Lagrangian. (Terpstra's results imply that such
quadratic forms exist only if $\ell\geq 3$ and $d\geq 3$). This indicates that, in general, the sharp bounds 
we obtain cannot be obtained using null-Lagrangians.

We allow for more general fields $\BE(\Bx)$, namely those with $m$ components $E_r(\Bx),~r=1,2,\ldots{},m$,  that derive from some real or complex
potential $\BU(\Bx)$, with $\ell$ components $U_1(\Bx),\ldots{},U_{\ell}(\Bx)$ through the equations
\beq E_{r}(\Bx)= \sum_{q=1}^{\ell}L_{rq}U_{q}(\Bx), \eeq{1.1}
for $r=1,2,\ldots{},m$, where $L_{rq}$ is the differential operator
\beq L_{rq}=A_{rq}+\sum_{h=1}^t\sum_{a_1,\ldots{},a_h=1}^{d} A_{rqh}^{a_1\ldots{}a_h}
{\Md\over\Md x_{a_1}}{\Md\over\Md x_{a_2}}\ldots{\Md\over\Md x_{a_h}}
\eeq{1.2}
of order $t$ in a space of dimension $d$ with real or complex valued constant coefficients $A_{rq}$ and
$A_{rqh}^{a_1\ldots{}a_h}$ , some of which may be zero. 

We will find that sharp inequalities on the integral \eq{0.5} can be obtained for certain boundary conditions when $f$ is a sharply $Q^*$-convex function that is quadratic.  A $Q^*$-convex function is defined as a function that satisfies 
\beq \lang f(\BE) \rang \geq f(\lang\BE\rang),
\eeq{1.2a}
for all periodic functions $\BE=\BL\BU$ that derive from a potential $\BU(\Bx)$ that is the sum of a polynomial $\BU^0(\Bx)$ and a periodic potential $\BU^1(\Bx)$ .
Here the angular brackets denote volume averages over the unit cell of periodicity. (By all periodic functions we also mean for all primitive unit cells of periodicity, including parallelepiped shaped ones). In the degenerate case, as pointed out to me by Marc Briane (private communication), there may be potentials that give rise to periodic fields, but are not expressible as the sum of a periodic part and a polynomial part, e.g. with $L=\Md/\Md x_1$
the potential $U=x_1\cos(x_2)$ is not so expressible. However such potentials are not of interest to us. 

The function $f$ is sharply  $Q^*$-convex  if in addition one has the equality
\beq \lang f(\underline{\BE}) \rang = f(\lang\underline{\BE}\rang),
\eeq{1.2b}
for some non-constant periodic function $\underline{\BE}=\BL\underline{\BU}$ that derives from a potential $\underline\BU(\Bx)$ that is the sum of a polynomial $\underline{\BU}^0(\Bx)$ and a periodic potential $\underline{\BU}^1(\Bx)$ . We call the fields $\underline{\BE}$ and its potential $\underline\BU(\Bx)$ 
$Q^*$-special fields.  When the elements $A_{rq}$ vanish and the coefficients  $A_{rqh}^{a_1\ldots{}a_h}$ are zero for all $h\ne t$
then $Q^*$-convexity may be equivalent to quasiconvexity, but otherwise they are not equivalent: see \cite{Ball:1981:NLW} for a precise definition of quasiconvexity.

Note that there could be a sequence of functions $\BE_1(\Bx), \BE_2(\Bx), \BE_3(\Bx),\ldots, \BE_n(\Bx)\ldots$  each having the same average value $\BE^0$, but not converging to $\BE^0$,
such that 
\beq \lim_{n\to\infty} \lang f(\BE_n) \rang = f(\BE^0)
\eeq{1.2c}
in which case we would call $f$ marginally $Q^*$-convex, but not sharply $Q^*$-convex unless there existed a field $\underline{\BE}$ such that the equality \eq{1.2b} held. In this paper we are not
interested in marginally $Q^*$-convex functions which are not sharply $Q^*$-convex.

In a nutshell, the main argument presented in this paper can be summarized as follows. When $\BE$ is the $Q^*$-special field $\underline{\BE}=\BL\underline{\BU}$ then we can
use integration by parts to show that the integral \eq{0.5} can be evaluated exactly in terms of boundary values, for all compact regions $\GO$ and in particular for regions $\GO$ within a cell $C$ of periodicity. Then we modify $\BU$ within $\GO$ while keeping $\BE=\underline{\BE}$ in $C\setminus\GO$, and maintaining the boundary
conditions so that  $\BE=\BL\BU$ holds weakly, including across the boundary $\Md\GO$. The field $\BE$ is extended outside $C$ to be periodic with unit cell $C$. Then the inequality \eq{1.2a} must
hold. This inequality (and an additional supplementary condition which ensures that $f(\lang\BE\rang)=f(\lang\underline{\BE}\rang)$) shows the integral of $f(\BE)$ over $\GO$ can only increase when $\BE$ is modified in this way.  This gives the desired sharp inequality on the integral. 

An associated problem, which we are also interested in, is to obtain sharp bounds on integrals of the form
\beq  \int_\GO g(\BJ(\Bx))\,d{\Bx}\ \eeq{1.3}
for fields $\BJ(\Bx)$ with $m$ real or complex components $J_r(\Bx),~r=1,2,\ldots{},m$,  that satisfy the differential constraints 
\beq \sum_{r=1}^{m}L_{qr}^\dag J_{r}(\Bx)=0
\eeq{1.4}
for $q=1,2,\ldots{},\ell$, where
\beq L_{qr}^\dag=\overline{A_{rq}}+\sum_{h=1}^t\sum_{a_1,\ldots{},a_h=1}^{d} (-1)^h \overline{A_{rqh}^{a_1\ldots{}a_h}}
{\Md\over\Md x_{a_1}}{\Md\over\Md x_{a_2}}\ldots{\Md\over\Md x_{a_h}},
\eeq{1.4a}
in which the bar denotes complex conjugation. Observe that the operators in $\eq{1.2}$ and $\eq{1.4a}$ are formal adjoints. 

Again, using essentially the same argument, we will find that sharp inequalities on the integral \eq{1.3} can be obtained for certain boundary conditions when $g$ is a 
sharply $Q^*$-convex function that is quadratic. In this setting a 
$Q^*$-convex function $g$ is defined as a function that satisfies 
\beq \lang g(\BJ) \rang \geq g(\lang\BJ\rang),
\eeq{1.4b}
for all periodic functions $\BJ$ satisfying $\BL^\dag\BJ=0$.  The function is sharply $Q^*$-convex  if in addition one has the equality
\beq \lang g(\dund{\BJ}) \rang = g(\lang\dund{\BJ}\rang),
\eeq{1.4c}
for some non-constant periodic function $\dund{\BJ}(\Bx)$ satisfying $\BL^\dag\dund{\BJ}=0$. We will call $\dund{\BJ}(\Bx)$
a $Q^*$-special field.

Presumably, the assumption that $f$ and $g$ are quadratic is not essential, but if they are not it becomes a difficult task to find ones which are sharply $Q^*$-convex.
Also there should be some generalization of the theory presented here to allow for $Q^*$-special fields  $\underline{\BE}$  and $\dund{\BJ}(\Bx)$ which are not
periodic, but we do not investigate this here. 

The theory developed here generalizes an example given in \cite{Kang:2013:BVF3d}, reviewed in the next section, that itself stemmed from developments in the calculus of variations, the theory
of topology optimization, and the theory of composites: see the books \cite{Cherkaev:2000:VMS, Allaire:2002:SOH, Milton:2002:TOC, Tartar:2009:GTH}. 
A key component of this is the Fourier space methods developed by Tartar and Murat \cite{Tartar:1979:CCA, Murat:1985:CVH, Tartar:1985:EFC} in their theory of compensated compactness for determining the quasiconvexity of quadratic forms.

\section{An example}
\setcounter{equation}{0}
To make the general analysis easier to follow we first review an example given in \cite{Kang:2013:BVF3d}, which was used to obtain sharp estimates of the volume occupied by an inclusion in a body from electrical impedance tomography measurements made at the surface of the body. This example serves to introduce the central arguments and the notion of  $Q^*$-special fields, and reviews the method of  Tartar and Murat  \cite{Tartar:1979:CCA, Murat:1985:CVH, Tartar:1985:EFC} for establishing the quasiconvexity of quadratic forms.

Let us consider in dimension $d=3$ a $3\times 3$ real valued matrix valued field $\BJ(\Bx)$ satisfying $\Div\BJ=0$. (Thus following the usual convention, but opposite
to the convention adopted in  \cite{Kang:2013:BVF3d,Milton:2002:TOC} its three rows, not columns, are each divergence free). Alternatively $\BJ(\Bx)$ could be regarded as a $9$ component vector which corresponds to the case $m=9$ and $t=1$ in \eq{1.4} and \eq{1.4a}.
We want to show, that for certain fluxes $\Bq=\BJ\Bn$ at the boundary $\Md\GO$
one can obtain sharp lower bounds on the integral \eq{1.3} when
\beq g(\BJ)=\Tr(\BJ^2)+\Tr(\BJ^T\BJ)-[\Tr(\BJ)]^2. \eeq{ex.1}
This special function, introduced in section 25.7 of \cite{Milton:2002:TOC}, is known to be quasiconvex (which in this case is equivalent to $Q^*$ convexity), in the sense that the inequality
\beq \lang g(\BJ)\rang-g (\lang\BJ\rang)\geq 0
\eeq{ex.2}
holds for all periodic functions $\BJ(\Bx)$ satisfying $\Div\BJ=0$ (as we are dealing with divergence free fields rather than gradients this is the appropriate definition of
quasiconvexity) where the angular brackets denote volume averages over the unit cell of periodicity. Following the 
ideas of Tartar and Murat  \cite{Tartar:1979:CCA, Murat:1985:CVH, Tartar:1985:EFC} the reason 
can be seen in Fourier space where the inequality (by Parseval's theorem, since the function $g$ is real and quadratic) takes the form
\beq  \sum_{\Bk\ne 0}g(\Real[\widehat{\BJ}(\Bk)])+g(\Imag[\widehat{\BJ}(\Bk)])\geq 0,
\eeq{ex.3}
which holds if  $g(\BH)$ is non-negative for all rank two real matrices $\BH$, and in particular for the matrices $\Real[\widehat{\BJ}(\Bk)]$ and $\Imag[\widehat{\BJ}(\Bk)]$
which are at most rank two since $\Div\BJ=0$ implies $\widehat{\BJ}(\Bk)\Bk=0$. To prove that $g(\BH)$ is non-negative for all rank two real matrices (i.e. rank 2 convex), one first notes that $g(\BH)$ is rotationally invariant in the sense that
\beq g(\BR^T\BH\BR)=g(\BH) 
\eeq{ex.4}
holds for all rotations $\BR$. This implies it suffices to check non-negativity with matrices $\BH$ such that $\BH[0,0,1]^T=0$, i.e. of the form
\beq \BH=\pmatrix{h_{11} & h_{12} & 0 \cr
h_{21} & h_{22} &  0 \cr
h_{31} & h_{32} & 0},
\eeq{ex.5}
 in which case
\beq g(\BH)=(h_{11}-h_{22})^2+(h_{12}+h_{21})^2+h_{31}^2+h_{32}^2
\eeq{ex.6}
is clearly non-negative and zero only when $h_{11}=h_{22}$, $h_{12}=-h_{21}$, and $h_{31}=h_{32}=0$. More generally, $g(\BH)$ will be zero if and only if
the rank two matrix $\BH$ has matrix elements of the form
\beq  H_{ij}=\Ga_0(k_ik_j-\Gd_{ij}\Bk\cdot\Bk)+\Gb_0\Ge_{ijm}k_m, 
\eeq{ex.7}
for some constants $\Ga_0$ and $\Gb_0$, and for some vector $\Bk$ which is a vector such that $\BH\Bk=0$, where $\Ge_{ijm}$ is the completely
antisymmetric Levi-Civita tensor taking the value $+1$ when $ijm$ is an even permutation of $123$, $-1$ when it is an odd permutation, and $0$ otherwise. 

 This establishes \eq{ex.2} and shows one has the equality
\beq  \lang g(\dund{\BJ})\rang-g (\lang\dund{\BJ}\rang)= 0,
\eeq{ex.8}
for  $Q^*$-special fields $\dund{\BJ}(\Bx)$ of the form
\beq \smash{\dund{J}}_{\ell k}=J^0_{\ell k}+\frac{\Md^2\Ga}{\Md x_\ell \Md x_k}-\Gd_{\ell k}\Delta\Ga+\Ge_{\ell km}\frac{\Md\Gb}{\Md x_m},
\eeq{ex.9}
where the $J^0_{\ell k}$ are elements of a constant matrix $\BJ^0$ and $\Ga(\Bx)$ and $\Gb(\Bx)$ are arbitrary periodic functions (in this case $\dund{\BJ}(\Bx)$ has Fourier components of the form \eq{ex.7}). This can be established directly by noting that
\beq g(\dund{\BJ})=\Tr(\dund{\BJ}^T\dund{\BE}),
\eeq{ex.10}
where
\beq \dund{\BE}=\dund{\BJ}+\dund{\BJ}^T-\Tr(\dund{\BJ})\BI=\BE^0+2\Grad\Grad\Ga, \quad \BE^0=\BJ^0+(\BJ^0)^T-\Tr(\BJ^0)\BI,
\eeq{ex.11}
and so \eq{ex.8} follows by integration by parts:
\beqa \lang g(\dund{\BJ})\rang=\lang \Tr(\dund{\BJ}^T\dund{\BE})\rang & = & \Tr\{(\BJ^0)^T[\BJ^0+(\BJ^0)^T-\Tr(\BJ^0)\BI]\}+2\lang\Tr[(\dund{\BJ})^T\Grad\Grad\Ga]\rang \nonum
& = &  g(\BJ^0)-2\lang(\Div\dund{\BJ})^T\Grad\Ga]\rang=g(\lang\dund{\BJ}\rang).
\eeqa{ex.12}

Now suppose $\BJ(\Bx)$ defined within $\GO$  satisfies the boundary condition 
$\Bq=\BJ\Bn$ at the boundary $\Md\GO$ where $q_\ell$ has components
\beq  q_\ell=\left(J^0_{\ell k}+\frac{\Md^2\Ga}{\Md x_\ell \Md x_k}-\Gd_{\ell k}\Delta\Ga+\Ge_{\ell k m}\frac{\Md\Gb}{\Md x_m}\right)n_k,
\eeq{ex.13}
for some functions $\Ga(\Bx)$ and $\Gb(\Bx)$ defined in the neighborhood of $\Md\GO$,  and for some constants $J^0_{\ell k}$. The key point is that one can extend  $\Ga(\Bx)$ and $\Gb(\Bx)$ beyond this boundary, so that they are periodic (with a unit cell of periodicity $C$ containing $\GO$) and define a field $\dund{\BJ}(\Bx)$ given by \eq{ex.9}.
The extension of $\BJ(\Bx)$ equals $\dund{\BJ}(\Bx)$ on $C\setminus\GO$ and is defined periodically over the entire domain
with periodic cell $C$.  The boundary conditions
\eq{ex.13} ensure that this periodic function satisfies $\Div\BJ(\Bx)$ in a weak sense and therefore \eq{ex.2} implies
\beq \int_{C\setminus\GO}g(\dund{\BJ}(\Bx))\,d\Bx+\int_{\GO}g(\BJ(\Bx))\,d\Bx\geq |C|g(\lang\BJ\rang),
\eeq{ex.14}
where $|C|$ is the volume of $C$. Also \eq{ex.8} implies 
\beq \int_{C\setminus\GO}g(\dund{\BJ}(\Bx))\,d\Bx+\int_{\GO}g(\dund{\BJ}(\Bx))\,d\Bx=|C|g(\lang\dund{\BJ}\rang).
\eeq{ex.15}
Furthermore since $\Div\BJ=0$ weakly, integration by parts implies that 
 \beq \int_{C}J_{k\ell}\,d\Bx=\int_{C}\sum_{i=1}^3\frac{\Md x_k}{\Md x_i}J_{i\ell}\,d\Bx
=\int_{\Md C} \sum_{i=1}^3 x_k n_i \smash{\dund{J}}_{i\ell}\,dS=\int_{C}\smash{\dund{J}}_{k\ell}\,d\Bx,
\eeq{ex.16}
and thus it follows that $\lang\BJ\rang=\lang\dund{\BJ}\rang$. So subtracting \eq{ex.15} from \eq{ex.14} gives 
\beq \int_{\GO}g(\BJ(\Bx))\,d\Bx\geq \int_{\GO}g(\dund{\BJ}(\Bx))\,d\Bx.
\eeq{ex.17}
Using \eq{ex.10} and integrating by parts allows us to evaluate the right hand side:
\beqa \int_{\GO}g(\dund{\BJ}(\Bx))\,d\Bx & = & \int_{\GO}\Tr[(\dund{\BJ})^T\Grad(\Bx^T\BE_0+2\Grad\Ga)]\,d\Bx=
\int_{\Md\GO}\Bn\cdot(\dund{\BJ})^T(\BE_0\Bx+2\Grad\Ga)\,dS \nonum 
& = & \int_{\Md\GO}\Bq\cdot(\BE_0\Bx+2\Grad\Ga)\,dS,
\eeqa{ex.18}
where $\BE_0$ is given by \eq{ex.11}.

In summary, for boundary fluxes of the form \eq{ex.13} we obtain the inequality
\beq \int_{\GO}g(\BJ(\Bx))\,d\Bx\geq \int_{\Md\GO}\Bq\cdot[(\BJ^0+(\BJ^0)^T-\Tr(\BJ^0)\BI)\Bx+2\Grad\Ga]\,dS,
\eeq{ex.19}
which is sharp, it being satisfied as an equality when $\BJ(\Bx)=\dund{\BJ}(\Bx)$ within $\GO$. In fact there is a huge range of fields for which one has equality since
one is free to change $\Ga(\Bx)$ and $\Gb(\Bx)$ in the interior of $\GO$ so long as the right side of \eq{ex.13} remains unchanged at the boundary $\Md\GO$.

\section{Conditions for $Q^*$-convexity}
\setcounter{equation}{0}
Here we solve the problem of characterizing those quadratic functions $f$ or $g$ which are $Q^*$-convex. This is a straightforward extension
of the ideas developed
by Tartar and Murat \cite{Tartar:1979:CCA, Murat:1985:CVH, Tartar:1985:EFC} in their theory of compensated compactness for characterizing quadratic quasiconvex functions. 
Again the key step is to study the inequality defining  quadratic $Q^*$-convex functions in the Fourier domain, where the differential constraints on the fields become algebraic constraints,  
and to use Parseval's theorem.

Since $f$ and $g$ are quadratic we can let
\beqa f(\BE(\Bx)) & = & \BE\cdot\BS\BE\equiv \sum_{r=1}^m\sum_{s=1}^m \overline{E}_r(\Bx)S_{rs}E_s(\Bx), \nonum
g(\BJ(\Bx)) & = & \BJ\cdot\BT\BJ\equiv \sum_{r=1}^m\sum_{s=1}^m \overline{J}_r(\Bx)T_{rs}J_s(\Bx).
\eeqa{1.5}
where the $S_{rs}$ or $T_{rs}$ are the real or complex valued elements of some Hermitian $m\times m$ matrix $\BS$ or $\BT$.

We first consider periodic functions $\BE(\Bx)$  and $\BJ(\Bx)$ and potentials $\BU(\Bx)$ that can be expressed in the form
\beq \BU(\Bx)=\BU^0(\Bx)+\BU^1(\Bx), \eeq{1.5a}
where $\BU^1(\Bx)$ is periodic with zero average value and $\BU^0(\Bx)$ is a polynomial with elements
\beq U^0_{q}=B_q+\sum_{h=1}^t\sum_{a_1,\ldots{},a_h=1}^{d} B_{qh}^{a_1\ldots{}a_h}
{x_{a_1}}{x_{a_2}}\ldots{x_{a_h}}, \eeq{1.5aa}
where the coefficients $B_q$ and $B_{qh}^{a_1\ldots{}a_h}$ are chosen so that 
\beq E^0_{r}=\sum_{q=1}^{\ell}L_{rq}U^0_{q}(\Bx) \eeq{1.5ab}
is a constant for $r=1,2,\ldots{},m$. We let $\CE^0$ denote the vector space spanned by all constant fields $\BE^0$  with elements $E^0_r$ expressible in the form \eq{1.5ab}.

We expand $\BE(\Bx)$, $\BU^1(\Bx)$ and $\BJ(\Bx)$  in a Fourier series
\beqa \BE(\Bx) & = & \lang\BE\rang+\sum_{\Bk\ne 0}\widehat{\BE}(\Bk)e^{i\Bk\cdot\Bx} \nonum
 \BU^1(\Bx) & = & \sum_{\Bk\ne 0}\widehat{\BU^1}(\Bk)e^{i\Bk\cdot\Bx}  \nonum
 \BJ(\Bx) & = & \lang\BJ\rang+\sum_{\Bk\ne 0}\widehat{\BJ}(\Bk)e^{i\Bk\cdot\Bx}, 
\eeqa{1.5b}
where the sum is over all $\Bk$ in the reciprocal lattice space (the reciprocal lattice consists of all $\Bk$ where the Fourier transform of the periodic functions $\BE(\Bx)$, $\BU^1(\Bx)$ and $\BJ(\Bx)$ have their natural support so that their primitive unit cell of periodicity is also a cell of periodicity of the wave
$e^{i\Bk\cdot\Bx}$, i.e. $e^{i\Bk\cdot\Ba}=1$ for all primitive lattice vectors $\Ba$ of these periodic functions).
From these expansions we see that the differential constraints \eq{1.1} or \eq{1.4} imply that for all $\Bk\ne 0$ in the reciprocal lattice
\beq \widehat{E}_{r}(\Bk)=\sum_{q=1}^{\ell}\widehat{L}_{rq}(\Bk)\widehat{U^1_{q}}(\Bk),\quad\sum_{r=1}^{m}\widehat{L}_{qr}^\dag(\Bk)\widehat{J}_{r}(\Bk)=0,
\eeq{1.6}
where $\widehat{L}_{qr}(\Bk)$ and $\widehat{L}_{qr}^\dag(\Bk)$ are the matrix elements
\beqa 
\widehat{L}_{rq}(\Bk) & = & A_{rq}+\sum_{h=1}^t\sum_{a_1,\ldots{},a_h=1}^{d} i^hA_{rqh}^{a_1\ldots{}a_h}
k_{a_1}k_{a_2}\ldots{}k_{a_h}, \nonum
\widehat{L}_{qr}^\dag (\Bk) & = & \overline{A_{rq}}+\sum_{h=1}^t\sum_{a_1,\ldots{},a_h=1}^{d} (-i)^h\overline{A_{rqh}^{a_1\ldots{}a_h}}
k_{a_1}k_{a_2}\ldots{}k_{a_h}
\eeqa{1.7}
of the $m\times m$ matrix $\widehat{\BL}(\Bk)$ and its adjoint $\widehat{\BL}^\dag (\Bk)$, defined for all $\Bk\in\RR^d$. In other words, the differential 
constraints imply that $\widehat{\BE}(\Bk)$ is in the range $\CE_{\Bk}$ of $\widehat{\BL}(\Bk)$ and that $\widehat{\BJ}(\Bk)$ is in the null-space $\CJ_{\Bk}$ of $\widehat{\BL}^\dag(\Bk)$,
where again these spaces are defined  for all $\Bk\in\RR^d$ and not just those $\Bk$ in the reciprocal lattice. We let $\BGG_1(\Bk)$
denote the projection onto $\CE_{\Bk}$ and $\BGG_2(\Bk)$ denote the projection onto $\CJ_{\Bk}$. Since these spaces are orthogonal complements it follows that
\beq \BGG_1(\Bk)+\BGG_2(\Bk)=\BI,\quad \BGG_1(\Bk)\BGG_2(\Bk)=\BGG_2(\Bk)\BGG_1(\Bk)=0,\quad\forall \Bk\in\RR^d
\eeq{1.7a}
where $\BI$ is the $m\times m$ identity matrix.

By substituting the Fourier expansion \eq{1.5b} for $\BE(\Bx)$ in \eq{1.2a} we see (by Parseval's theorem) that \eq{1.2a} holds if and only if the expression
 \beq  \lang\BE\cdot\BS\BE\rang  -  \lang\BE\rang\cdot\BS\lang\BE\rang
= \sum_{\Bk\ne 0}\widehat{\BE}(\Bk)\cdot\BS\widehat{\BE}(\Bk)= \sum_{\Bk\ne 0}f(\widehat{\BE}(\Bk))
 \nonum
\eeq{1.9}
is non-negative. A necessary and sufficient condition for this to hold for all possible primitive unit cells is that for all non-zero $\Bk\in\RR^d$ 
\beq f(\BH)\geq 0\quad{\rm for\,\,all}\,\,\BH\in\CE_{\Bk},
\eeq{1.10}
or equivalently that
\beq \BGG_1(\Bk)\BS\BGG_1(\Bk)\geq 0, \quad \forall \Bk\in\RR^d,\,\, \Bk\ne 0,
\eeq{1.10a} 
where the inequality holds in the sense of quadratic forms. This is an algebraic condition that can be checked numerically, and in some cases analytically.  

If $\BS$ is real (and hence symmetric),
then \eq{1.9} reduces to
\beq  \lang\BE\cdot\BS\BE\rang  -  \lang\BE\rang\cdot\BS\lang\BE\rang
=\sum_{\Bk\ne 0}\Real[\widehat{\BE}(\Bk)]\cdot\BS\Real[\widehat{\BE}(\Bk)]
+\sum_{\Bk\ne 0}\Imag[\widehat{\BE}(\Bk)]\cdot\BS\Imag[\widehat{\BE}(\Bk)].
\eeq{1.11}
If in addition the coefficients $A_{rq}$ are real and the coefficients
$A_{rqh}^{a_1\ldots{}a_h}$ are real when $h$ is even and purely imaginary when $h$ is odd
then $\widehat{\BL}(\Bk)$ is real and $\BH\in\CE_{\Bk}$ if and only if the real and imaginary parts of 
$\BH$ lie in $\CE_{\Bk}$. Thus, in this case, to guarantee \eq{1.2a} it suffices that \eq{1.10} holds for real $\BH$,  for all $\Bk\ne 0$.

Similarly, we look for Hermitian matrices $\BT$ which are $Q^*$-convex in the sense that \eq{1.4b} holds,
and a necessary and sufficient condition for this is that 
\beq g(\BH)\geq 0\quad{\rm for\,\,all}\,\,\BH\in\CJ_{\Bk},
\eeq{1.13}
or equivalently that
\beq \BGG_2(\Bk)\BT\BGG_2(\Bk)\geq 0, \quad \forall \Bk\in\RR^d,\,\,  \Bk\ne 0, \eeq{1.14}
where again the inequality hold in the sense of quadratic forms. 
If $\BT$ is real and  the coefficients $A_{rq}$ are real and the coefficients
$A_{rqh}^{a_1\ldots{}a_h}$ are real when $h$ is even and purely imaginary when $h$ is odd, then it suffices to check
\eq{1.13} holds for real $\BH$. 

\section{Sharply $Q^*$-convex quadratic functions and their associated $Q^*$-special fields}
\setcounter{equation}{0}
Here we are interested in sharply $Q^*$-convex functions $f(\BE)$ for which one has the equality \eq{1.2b}
for some non-constant periodic function $\underline{\BE}=\BL\underline{\BU}$ that derives from a potential $\underline\BU(\Bx)$ that is the sum of a polynomial $\underline{\BU}^0(\Bx)$ and a periodic potential $\underline{\BU}^1(\Bx)$. We will see that associated with $\underline{\BE}$ is a companion field $\underline{\BJ}=\BS\underline{\BE}$ satisfying $\BL^\dag\underline{\BJ}=0$. As remarked in the next section this allows us to express the integral of $f(\underline{\BE})=\underline{\BE}\cdot\underline{\BJ}$ over $\GO$ in terms of boundary values. Similarly we are  interested in sharply $Q^*$-convex functions $g(\BJ)$ for which one has the equality \eq{1.4c} for some non-constant periodic function $\dund{\BJ}$  satisfying $\BL^\dag\dund{\BJ}=0$.  We will see that associated with $\dund{\BJ}$ is a companion field $\dund{\BE}$ and a companion potential $\dund{\BU}$ such that $\dund{\BE}=\BT\dund{\BJ}=\BL\dund{\BU}$. This will allow us to express the integral of $g(\dund{\BJ})=\dund{\BJ}\cdot\dund{\BE}$ over $\GO$ in terms of boundary values.

The $Q^*$-convexity conditions \eq{1.10a} and \eq{1.14} are clearly satisfied when $\BS$ and $\BT$ are positive definite, but in this case one only
has equality in \eq{1.10} and \eq{1.13} when $\BH=0$.  
For any $\Bk\ne 0$ let $\CS_{\Bk}$ denote the subspace of all $m$ dimensional vectors $\BH\in\CE_{\Bk}$ for which 
one has equality in \eq{1.10}, and let $\CU_{\Bk}$ denote the subspace of all $\ell$ dimensional vectors $\BG$ such that $\widehat{\BL}(\Bk)\BG\in\CS_{\Bk}$.
We call $f$ sharply $Q^*$-convex if, for some $\Bk\ne 0$, $\CS_{\Bk}$ contains at least one non-zero vector $\BH$. Note
that if $\BH\in\CS_{\Bk}$ then 
\beq \BH\cdot\BGG_1\BS\BGG_1\BH=0, \eeq{2.1}
and since from \eq{1.10a} $\BGG_1\BS\BGG_1$ is positive semi-definite we deduce that $\BH$ must in fact be a null-vector of this matrix, which implies that
\beq \BGG_1\BS\BH=0, \quad{i.e.}\,\,\BS\BH\in \CJ_{\Bk}, \eeq{2.2}
and hence that
\beq \BS\widehat{\BL}(\Bk)\BG\in \CJ_{\Bk},\quad\forall\,\BG\in \CU_{\Bk}.
\eeq{2.3}
Associated with a sharply $Q^*$-convex $f$ are periodic potentials $\underline{\BU}^1(\Bx)$ expressible in the form
\beq \underline{\BU}^1(\Bx)  = \sum_{\Bk\ne 0}\widehat{\underline{\BU}^1}(\Bk)e^{i\Bk\cdot\Bx},
\eeq{2.4}
where $\widehat{\underline{\BU}^1}(\Bk)\in \CU_{\Bk}$ for all $\Bk\ne 0$.  The unit cell of periodicity has to be
chosen so the $\Bk$ in the reciprocal lattice include some $\Bk$ such that  $\CS_{\Bk}$ contains at least one non-zero vector.

Let us introduce the companion $Q^*$-special fields
\beqa \underline{\BE}(\Bx) & = & \underline{\BE}^0+\sum_{\Bk\ne 0}\widehat{\underline{\BE}}(\Bk)e^{i\Bk\cdot\Bx}, \nonum
\underline{\BJ}(\Bx) & = & \BS\underline{\BE}^0+\sum_{\Bk\ne 0}\widehat{\underline{\BJ}}(\Bk)e^{i\Bk\cdot\Bx}, 
\eeqa{2.5}
in which $\underline{\BE^0}\in \CE^0$ so that $\underline{\BE^0}=\BL\underline{\BU}^0(\Bx)$ for some polynomial potential $\underline{\BU}^0(\Bx)$ and where
\beq \widehat{\underline{\BE}}(\Bk)=\widehat{\BL}(\Bk)\widehat{\underline{\BU}^1}(\Bk), \quad 
 \widehat{\underline{\BJ}}(\Bk)= \BS\widehat{\underline{\BE}}(\Bk).
\eeq{2.6}
With $ \underline{\BU}(\Bx)=\underline{\BU}^0(\Bx)+ \underline{\BU}^1(\Bx)$ these $Q^*$-special fields satisfy
\beq \underline{\BJ}(\Bx)=\BS\underline{\BE}(\Bx),\quad \underline{\BE}=\BL\underline{\BU},\quad \BL^\dag\underline{\BJ}=0,
\eeq{2.7}
and
\beq  \lang\underline{\BE}\cdot\BS\underline{\BE}\rang  -  \lang\underline{\BE}\rang\cdot\BS\lang\underline{\BE}\rang
= \sum_{\Bk\ne 0}\widehat{\underline{\BE}}(\Bk)\cdot\BS\widehat{\underline{\BE}}(\Bk)=0,
\eeq{2.7a}
so that \eq{1.2b} is satisfied.

Similarly for all $\Bk\ne 0$ let $\CT_{\Bk}$ denote the subspace of all $m$ dimensional vectors $\BH\in\CJ_{\Bk}$ for which 
one has equality in \eq{1.13}. We call $g$ sharply $Q^*$-convex if, for some $\Bk\ne 0$, $\CT_{\Bk}$ contains at least one non-zero vector $\BH$.
Analogous to \eq{2.2}, \eq{1.14} implies
\beq \BGG_2\BT\BH=0, \quad{i.e.}\,\,\BT\BH\in \CE_{\Bk}. \eeq{2.8}
Since $\CE_{\Bk}$ is the range of $\widehat{\BL}(\Bk)$, there exists a (possibly non-unique) $\ell$ dimensional vector $\BG$ such that
\beq \BT\BH=\widehat{\BL}(\Bk)\BG. \eeq{2.9}
Now associated with a sharply $Q^*$-convex $g$ are the $Q^*$-special fields
\beqa 
\dund{\BJ}(\Bx) & = & \dund{\BJ}^0+\sum_{\Bk\ne 0}\widehat{\dund{\BJ}}(\Bk)e^{i\Bk\cdot\Bx}, \nonum
\dund{\BE}(\Bx) & = & \BT\dund{\BJ}^0+\sum_{\Bk\ne 0}\widehat{\dund{\BE}}(\Bk)e^{i\Bk\cdot\Bx}, \nonum  
\eeqa{2.10}
where 
\beq \widehat{\dund{\BJ}}(\Bk)\in\CT_{\Bk},\quad
\widehat{\dund{\BE}}(\Bk)=\BT\widehat{\dund{\BJ}}(\Bk)\in \CE_{\Bk},
\eeq{2.11}
and $\BT\dund{\BJ}^0\in\CE^0$ so that $\BT\dund{\BJ^0}=\BL\dund{\BU}^0(\Bx)$ for some polynomial potential $\dund{\BU}^0(\Bx)$. The unit cell of periodicity has to be
chosen so the $\Bk$ in the reciprocal lattice include some $\Bk$ such that  $\CT_{\Bk}$ contains at least one non-zero vector. We choose the 
periodic potential  
\beq \dund{\BU}^1(\Bx)  = \sum_{\Bk\ne 0}\widehat{\dund{\BU}^1}(\Bk)e^{i\Bk\cdot\Bx},
\eeq{2.12}
so that its Fourier coefficients satisfy 
\beq  \widehat{\dund{\BE}}(\Bk)=\widehat{\BL}(\Bk)\widehat{\dund{\BU}^1}(\Bk),
\eeq{2.13}
for $\Bk\ne 0$. With $\dund{\BU}(\Bx)=\dund{\BU}^0(\Bx)+ \dund{\BU}^1(\Bx)$ these $Q^*$-special fields satisfy
\beq \dund{\BE}(\Bx)=\BT\dund{\BJ}(\Bx),\quad \dund{\BE}=\BL\dund{\BU},\quad \BL^\dag\dund{\BJ}=0,
\eeq{2.14}
and
\beq  \lang\dund{\BJ}\cdot\BT\dund{\BJ}\rang  -  \lang\dund{\BJ}\rang\cdot\BT\lang\dund{\BJ}\rang
= \sum_{\Bk\ne 0}\widehat{\dund{\BJ}}(\Bk)\cdot\BT\widehat{\dund{\BJ}}(\Bk)=0,
\eeq{2.15}
so that \eq{1.4c} is satisfied.

\section{Sharp inequalities on the integrals over $\GO$}
\setcounter{equation}{0}
Since the operators $\BL$ and $\BL^\dag$ are formal adjoints the quantities
\beqa f_0=\int_\GO f(\underline{\BE}(\Bx))\,d{\Bx}=\int_\GO \underline{\BE}(\Bx)\cdot\BS\underline{\BE}(\Bx)=\int_\GO \underline{\BE}(\Bx)\cdot\underline{\BJ}(\Bx),\nonum
 g_0=\int_\GO g(\dund{\BE}(\BJ))\,d{\Bx}=\int_\GO \dund{\BJ}(\Bx)\cdot\BT\dund{\BJ}(\Bx)=\int_\GO \dund{\BJ}(\Bx)\cdot\dund{\BE}(\Bx)
\eeqa{3.1}
can be computed in terms of boundary terms using integration by parts. Now we show that 
\beq \int_\GO f(\BE(\Bx))\,d{\Bx}\geq f_0,\quad \int_\GO g(\BJ(\Bx))\,d{\Bx}\geq g_0,
\eeq{3.2}
for all fields $\BE(\Bx)$ deriving from a potential $\BU(\Bx)$ that matches the appropriate boundary data of the potential $\underline{\BU}(\Bx)$ (generally involving
both $\BU(\Bx)$ and its derivatives when $t>1$)
or for all fields $\BJ(\Bx)$ that match the appropriate boundary data of
$\dund{\BJ}(\Bx)$, provided certain further supplementary conditions hold. First note that these inequalities are clearly sharp, being attained when $\BE(\Bx)=\underline{\BE}(\Bx)$ or when
 $\BJ(\Bx)=\dund{\BJ}(\Bx)$. 

To establish the first inequality in \eq{3.2} we first find a parallelepiped $C$ that contains $\GO$ and that is formed from an integer number of the primitive unit cells
of  $\underline{\BE}(\Bx)$. (If $\GO$ lies inside a primitive unit cell of  $\underline{\BE}(\Bx)$, then we can take $C$ as this primitive cell, but otherwise we need to join 
a set of these primitive cells together to obtain a parallelepiped that covers $\GO$). We extend $\BE(\Bx)$ outside $\GO$ so that 
it is periodic with $C$ as a unit cell. In this cell, but outside $\GO$, $\BE(\Bx)$ equals $\underline{\BE}(\Bx)$. The boundary data on  
the potential $\underline{\BU}(\Bx)$ are chosen so the equation $\underline{\BE}=\BL\underline{\BU}$ holds weakly across the boundary of $\GO$. (For example,
if $m=\ell=d=1$ and $L=\Md^2/\Md x_1^2$ then this would require continuity of both $U$ and $\Md U/\Md x_1$ at the interface.)
We extend the potential $\BU(\Bx)$ outside $C$ so 
that $\BU^1(\Bx)\equiv \BU(\Bx)-\underline{\BU}^0(\Bx)$ is $C$-periodic: if $\Bx_0$ is any lattice vector 
then for $\Bx\in C$ we set
\beq \BU(\Bx+\Bx_0)=\BU(\Bx)+\underline{\BU}^0(\Bx+\Bx_0)-\underline{\BU}^0(\Bx).
\eeq{3.3}
Defined in this way, the relation $\BE=\BL\BU$ holds in a weak sense, and
so the inequality \eq{1.2a} is satisfied, which we rewrite as
\beq  \int_\GO f(\BE)\,d\Bx+\int_{C\setminus\GO} f(\underline{\BE})\,d\Bx  \geq |C| f(\lang\BE\rang),
\eeq{3.4}
where $|C|$ is the volume of $C$. Also \eq{2.7a} holds, which we rewrite as
\beq   f_0+\int_{C\setminus\GO} f(\underline{\BE})\,d\Bx  = |C| f(\lang\underline{\BE}\rang),
\eeq{3.5}
so subtracting these equations gives
\beq  \int_\GO f(\BE)\,d\Bx\geq f_0+ |C|[f(\lang\BE\rang)- f(\lang\underline{\BE}\rang)].
\eeq{3.6}
This inequality in general requires us to know $\lang\BE\rang$ and $\lang\underline{\BE}\rang$. Given a $m$-dimensional constant vector $\BJ^0$, we have
\beq \int_C\BJ^0\cdot({\BE}-\underline{\BE})\,d\Bx=\int_C\BJ^0\cdot\BL(\BU-\underline{\BU})\,d\Bx=
\int_C\BJ^0\cdot\BA(\BU-\underline{\BU})\,d\Bx,
\eeq{3.7}
where $\BA$ is the matrix with elements $A_{rq}$. (In establishing the last equality in \eq{3.7} we have used integration by parts and the fact that
$\BU$ and $\underline{\BU}$ are equal in the vicinity of the boundary of $C$). Since this holds for all $\BJ^0$, we deduce that
\beq \lang\BE\rang= \lang\underline{\BE}\rang+\frac{1}{|C|}\int_C \BA(\BU-\underline{\BU})\,d\Bx.
\eeq{3.8}
Therefore a sufficient condition for $f(\lang\BE\rang)$ to be equal to $f(\lang\underline{\BE}\rang)$ is that
\beq \BS\BA=0, \eeq{3.9}
so that the range of $\BA$ is in the null space of $\BS$ (which if $\BS$ is non-singular requires that $\BA=0$). When this supplementary condition holds then
clearly \eq{3.6} implies the first inequality in \eq{3.2}

To establish the second inequality in \eq{3.2} we first find a parallelepiped $C$ that contains $\GO$ and that is formed from an integer number of the primitive unit cells
of  $\underline{\BJ}(\Bx)$. Then we extend $\BJ(\Bx)$ outside $\GO$ so that 
it is $C$-periodic, and within the
unit cell $C$ equals $\underline{\BJ}(\Bx)$ outside $\GO$. The boundary data on  
the field $\underline{\BJ}(\Bx)$ are chosen so the equation $\BL\underline{\BJ}=0$ holds weakly across the boundary of $\GO$. Then analogously to \eq{3.6} we have
\beq  \int_\GO g(\BJ)\,d\Bx\geq g_0+ |C|[g(\lang\BJ\rang)- g(\lang\underline{\BJ}\rang)].
\eeq{3.12}
This inequality in general requires us to know $\lang\BJ\rang$ and $\lang\underline{\BJ}\rang$. Given a $m$ dimensional constant vector $\BE^0\in \CE^0$ so that $\BE^0=\BL\BU^0(\Bx)$ for some polynomial potential $\BU^0(\Bx)$, we have
\beq \int_C ({\BJ}-\underline{\BJ})\cdot{\BE^0}\,d\Bx=\int_C ({\BJ}-\underline{\BJ})\cdot\BL{\BU^0}\,d\Bx=0,
\eeq{3.13}
where the last equality follows from integration by parts, using the fact that ${\BJ}=\underline{\BJ}$ in a vicinity of the boundary of $C$. 
Since this holds for all $\BE^0\in\CE^0$ we deduce that a sufficient condition for $g(\lang\BE\rang)$ to be equal to $g(\lang\underline{\BE}\rang)$ is that
the range of $\BT$ be a subset of $\CE^0$ and if this supplementary condition holds then
clearly \eq{3.12} implies the second inequality in \eq{3.2}.

\section{An algorithm for generating sharply $Q^*$-convex quadratic functions and their associated $Q^*$-special fields and for generating extremal quasiconvex functions}
\setcounter{equation}{0}
For applications one needs a way of generating $Q^*$-convex quadratic functions and their associated $Q^*$-special fields. This section shows how to do this, but is fairly technical
and so can be skipped be readers not interested in the details. However, it is important to emphasize that the approach presented here also solves, for the first time,
the problem of generating arbitrary  extremal quasiconvex functions.

To begin with let us suppose $\BS$ can be expressed in the form
\beq \BS=\BV-\Bv\Bv^\dag/\Ga,
\eeq{2.16}
where $\BV$ is Hermitian and positive definite, and $\Ga$ is real and positive. The $Q^*$-convexity of $f$ is equivalent to the inequality
\beq \BH\cdot\BV\BH\geq \overline{\Gg}(\BH\cdot\Bv)+\Gg(\overline{\BH\cdot\Bv})-\Ga|\Gg|^2
\eeq{2.17}
holding for all complex $\Gg$ and for all $\BH\in\CE_{\Bk}$, for all non-zero $\Bk\in\RR^d$, as can be established by taking
the maximum of the right hand side over $\Gg$ thereby recovering \eq{1.10}. Rewriting this inequality as 
\beq \Ga|\Gg|^2\geq \overline{\Gg}(\BH\cdot\Bv)+\Gg(\overline{\BH\cdot\Bv})-\BH\cdot\BV\BH,
\eeq{2.18}
and taking the maximum of the right hand side over $\BH\in\CE_{\Bk}$ we see that an equivalent condition is that
\beq \Ga\geq \Bv\cdot\BGG(\Bk)\Bv, \quad \forall \Bk\in\RR^d,\,\,\Bk\ne 0,
\eeq{2.19}
where, assuming the $\ell\times\ell$ matrix $\widehat{\BL}^\dag(\Bk)\BV\widehat{\BL}(\Bk)$ is non-singular for all $\Bk\ne 0$, 
\beq \BGG(\Bk)=\widehat{\BL}(\Bk)(\widehat{\BL}^\dag(\Bk)\BV\widehat{\BL}(\Bk))^{-1}\widehat{\BL}^\dag(\Bk).
\eeq{2.20}
 In the context of quasiconvexity, the condition \eq{2.19} was first derived in \cite{Milton:1990:CSP},
stimulated by a result of Kohn and Lipton \cite{Kohn:1988:OBE}.
The simpler derivation presented here was suggested by an anonymous referee of a later paper (see also section 24.9 in \cite{Milton:2002:TOC}).
If $\BV=\BI$ then $\BGG(\Bk)$ is equal to $\BGG_1(\Bk)$, the projection onto the range $\CE_{\Bk}$ of $\widehat{\BL}(\Bk)$. More generally $\BGG(\Bk)$
is defined to be the matrix such that for all vectors $\Bv$, 
\beq \BGG(\Bk)\Bv\in \CE_{\Bk}\quad{\rm and}\quad  \BV\BGG(\Bk)\Bv-\Bv \in \CJ_{\Bk},
\eeq{2.20a}
implying that
\beq \BGG(\Bk)=\BGG_1(\Bk)[\BGG_1(\Bk)\BV\BGG_1(\Bk)]^{-1}\BGG_1(\Bk),
\eeq{2.20ab}
where the inverse is to be taken on the space $\CE_{\Bk}$. (Thus $\BGG(\Bk)$ is the psuedoinverse of $\BGG_1(\Bk)\BV\BGG_1(\Bk)$). To see that \eq{2.20a} implies
\eq{2.20ab} we use the fact that $\CE_{\Bk}$ and $\CJ_{\Bk}$ are orthogonal subspaces. Hence \eq{2.20a} implies
\beq \BGG_1\BGG(\Bk)\Bv=\BGG(\Bk)\Bv, \quad\BGG_1[ \BV\BGG(\Bk)\Bv-\Bv]=0. \eeq{2.20ac}
Combining these gives
\beq \BGG_1\BV\BGG_1(\BGG(\Bk)\Bv)=\BGG_1\Bv,\eeq{2.20ad}
which implies \eq{2.20ab}, since it holds for all $\Bv$.

Clearly \eq{2.19} is just satisfied if
\beq \Ga=\sup_{\matrix{\Bk\in\RR^d \cr \Bk\ne 0}} \Bv\cdot\BGG(\Bk)\Bv,
\eeq{2.21}
and in this case $f$ is marginally $Q^*$-convex. However if this supremum is attained for some $\Bk\ne 0$ then $f$ is sharply $Q^*$-convex.  

As an example, consider fields of the form
\beq \BE=\pmatrix{\Grad U \cr U}, \eeq{2.21a}
where $U(\Bx)$ is a scalar potential. Suppose $\BV=\BI$ and 
\beq \Bv=\pmatrix{\Bt \cr 1}, \eeq{2.21b}
where $\Bt$ is a real $d$-component unit vector. Then $\BGG(\Bk)$ is the projection 
\beq \BGG(\Bk)=\BGG_1(\Bk)=\frac{1}{k^2+1}\pmatrix{\Bk \cr 1}\pmatrix{\Bk^T & 1}, \eeq{2.21c}
in which $k^2=\Bk\cdot\Bk$ and so 
\beq \Bv\cdot\BGG(\Bk)\Bv=\frac{(\Bt\cdot\Bk+1)^2}{k^2+1}
\eeq{2.21d}
takes its maximum value $\Ga=2$ when $\Bk=\Bt$. (It is clear that the maximum occurs when $\Bt$ and $\Bk$ are parallel and then it is just
a matter of simple algebra to show \eq{2.21d} is less than or equal to $2$.)
This example shows that the supremum can be attained, and in some cases may be attained only at one non-zero value of $\Bk$. 

Returning to the general case, let $\CK$ denote the set
of those $\Bk\ne 0$ which attain the supremum in \eq{2.21}. Then for $\Bk\in\CK$, $\Ga=\Bv\cdot\BGG(\Bk)\Bv$ and the set
$\CU_{\Bk}$ consists of those $\BG$ of the form
\beq \BG=a(\widehat{\BL}^\dag(\Bk)\BV\widehat{\BL}(\Bk))^{-1}\widehat{\BL}^\dag(\Bk)\Bv,
\eeq{2.22}
where $a$ is an arbitrary complex constant. Thus
\beq \BH\equiv \widehat{\BL}(\Bk)\BG=a\widehat{\BL}(\Bk)(\widehat{\BL}^\dag(\Bk)\BV\widehat{\BL}(\Bk))^{-1}\widehat{\BL}^\dag(\Bk)\Bv=a\BGG(\Bk)\Bv
\eeq{2.23}
satisfies
\beq \BS\BH=\BV\BH-a\Bv=-a\BGD(\Bk)\BV^{-1}\Bv,
\eeq{2.23a}
where
\beq \BGD(\Bk)\equiv\BV-\BV\BGG(\Bk)\BV 
\eeq{2.23b}
has the property that
\beq \widehat{\BL}^\dag(\Bk)\BGD(\Bk)=\widehat{\BL}^\dag(\Bk)\BV-\widehat{\BL}^\dag(\Bk)\BV\widehat{\BL}(\Bk)(\widehat{\BL}^\dag(\Bk)\BV\widehat{\BL}(\Bk))^{-1}\widehat{\BL}^\dag(\Bk)\BV=0,
\eeq{2.23c}
implying $\BS\BH\in\CJ_{\Bk}$. So the potential $\underline{\BU}^1(\Bx)$ and the $Q^*$-special fields  $\underline{\BE}(\Bx)$ and $\underline{\BJ}(\Bx)$ take the form
\beqa \underline{\BU}^1(\Bx) & =  & \sum_{\Bk\ne 0}a(\Bk)[(\widehat{\BL}^\dag(\Bk)\BV\widehat{\BL}(\Bk))^{-1}\widehat{\BL}^\dag(\Bk)]\Bv\,e^{i\Bk\cdot\Bx}, \nonum
\underline{\BE}(\Bx) & = & \underline{\BE}^0+\sum_{\Bk\in\ne 0}a(\Bk)\BGG(\Bk)\Bv\,e^{i\Bk\cdot\Bx}, \nonum
\underline{\BJ}(\Bx) & =  & \BS\underline{\BE}(\Bx) =\BS\underline{\BE}^0-\sum_{\Bk\ne 0}a(\Bk)\BGD(\Bk)\BV^{-1}\Bv\,e^{i\Bk\cdot\Bx},
\eeqa{2.25}
where $\underline{\BE}^0\in \CE^0$ and
\beq a(\Bk)=0\quad{\rm if}\,\,\Bk\notin\CK.
\eeq{2.25a}
The reciprocal lattice needs to be chosen so that it includes some $\Bk\in\CK$ to ensure the $a(\Bk)$ in \eq{2.25} are not all zero.

When $\BV=\BI$ then $\BGD(\Bk)$ is equal to the projection $\BGG_2(\Bk)$ onto the null-space $\CJ_{\Bk}$ of $\widehat{\BL}^\dag(\Bk)$. More generally $\BGD(\Bk)$
is defined to be the matrix such that for all vectors $\Bw$, 
\beq \BGD(\Bk)\Bw\in \CJ_{\Bk}\quad{\rm and}\quad  \BV^{-1}\BGD(\Bk)\Bw-\Bw \in \CE_{\Bk},
\eeq{2.25aa}
implying 
\beq \BGD(\Bk)=\BGG_2(\Bk)[\BGG_2(\Bk)\BV^{-1}\BGG_2(\Bk)]^{-1}\BGG_2(\Bk), \eeq{2.25ab}
where the inverse is to be taken on the space $\CJ_{\Bk}$. 

We now use the fact that for an appropriate choice of the scalar polynomial $q(\Bk)$,
\beq \BP(\Bk)=q(\Bk)(\widehat{\BL}^\dag(\Bk)\BV\widehat{\BL}(\Bk))^{-1}\widehat{\BL}^\dag(\Bk)\Bv
\eeq{2.25b}
is also polynomial in $\Bk$: to see this it suffices to take 
\beq q(\Bk)=\det[\widehat{\BL}^\dag(\Bk)\BV\widehat{\BL}(\Bk)],
\eeq{2.25ba}
although it is better to take the lowest degree polynomial which works. Then setting
\beq a(\Bk)=q(\Bk)\widehat{c}(\Bk),
\eeq{2.25c}
where the $\widehat{c}(\Bk)$ are the Fourier components of some scalar periodic function $c(\Bx)$, we see that
\beq \underline{\BU}^1(\Bx)=\BP(\Md/\Md x_1, \Md/\Md x_2,\ldots , \Md/\Md x_d)c(\Bx),
\eeq{2.25d}
which also gives the $Q^*$-special fields $\underline{\BE}=\underline{\BE}^0+\BL\underline{\BU}^1$ and $\underline{\BJ}=\BS\underline{\BE}$. The function
$c(\Bx)$ has to be chosen so that $\widehat{c}(\Bk)$ is zero if $\Bk\notin\CK$. (We do not have to worry about this latter constraint if it happens that
$\CK$ consists of all non-zero vectors in $\RR^d$).

We assumed $\BV$ was positive definite but in fact it is easy to check that the derivation goes through if the quadratic form associated with $\BV$ is $Q^*$-convex.
It could also be marginally $Q^*$-convex but in this case one needs to choose $\Bv$ so that 
\beq \widehat{\BL}^\dag(\Bk)\Bv\quad{\rm lies\,\,in\,\,the\,\,range\,\,of}\,\,\widehat{\BL}^\dag(\Bk)\BV\widehat{\BL}(\Bk),\quad \forall\Bk\ne 0,
\eeq{2.25e}
since otherwise the supremum in \eq{2.21} is surely infinite. Thus by induction, successively
setting
\beq \BV=\BV_0-\sum_{i=1}^{j-1}\Bv_i\Bv_i^\dag/\Ga_i, \quad\Bv=\Bv_{j},\quad \Ga=\Ga_{j},
\eeq{2.25ea}
for $j=1,2,\ldots, n$ (where $\BV_0$ has an associated $Q^*$-convex quadratic form, 
$\Bv=\Bv_{j}$ has to be chosen so that \eq{2.25e} is satisfied, and $\Ga=\Ga_{j}$ is given by \eq{2.21}) one finally reaches a point  
 such that the quadratic function $f$ associated with 
\beq \BS=\BV_0-\sum_{i=1}^n\Bv_i\Bv_i^\dag/\Ga_i, 
\eeq{2.25f}
is an extremal quadratic $Q^*$-convex function in the sense that one cannot subtract any non-zero positive semi-definite tensor from $\BS$ and retain its $Q^*$-convexity.
In this case, with $\BV=\BS$ the supremum in \eq{2.21} is infinite for all non-zero $\Bv\in\CC^m$. We emphasize that this also solves the outstanding problem of generating extremal quadratic quasiconvex functions (or equivalently "extremal translations"): see \cite{Milton:1990:CSP} and chapters 24 and 25 in \cite{Milton:2002:TOC}  for other insights into this problem. 

In the end one obtains $Q^*$-special potentials of the form
\beq  \underline{\BU}^1(\Bx)=\sum_{j=1}^n\BP_j(\Md/\Md x_1, \Md/\Md x_2,\ldots , \Md/\Md x_d)c_j(\Bx).
\eeq{2.25g}
Here the $n$ scalar periodic potentials $c_j(\Bx)$ have Fourier components satisfying 
\beq  \widehat{c}_j(\Bk)=0\quad{\rm if}\,\,\Bk\notin\CK_j,
\eeq{2.25h}
and, when one has made the substitutions \eq{2.25ea}, the $\BP_j$ are polynomials given by \eq{2.25b}  and the $\CK_j$ consist of those $\Bk$ which attain the supremum in \eq{2.21},

Similarly, supposing $\BT$ (the matrix defined by \eq{1.5} that is associated with the quadratic form $g$) can be expressed in the form
\beq \BT=\BV^{-1}-\Bw\Bw^\dag/\Gb,
\eeq{2.26}
where $\BV$ is positive definite, we see that $g$ is $Q^*$-convex if and only if
\beq \Gb\geq \Bw\cdot\BGD(\Bk)\Bw, \quad \forall \Bk\in\RR^d,\,\,\Bk\ne 0.
\eeq{2.27}
If we choose
\beq \Gb=\sup_{\matrix{\Bk\in\RR^d \cr \Bk\ne 0}}\Bw\cdot\BGD(\Bk)\Bw
\eeq{2.29}
then $g$ will be marginally $Q^*$-convex and sharply $Q^*$-convex  if this supremum is attained for some $\Bk\ne 0$.  Let $\CM$ denote the set
of those $\Bk\ne 0$ which attain the supremum. Then for $\Bk\in\CM$, $\Gb=\Bw\cdot\BGD(\Bk)\Bw$ and the set
$\CJ_{\Bk}$ consists of those $\BH$ of the form
\beq \BH=b\BGD(\Bk)\Bw,
\eeq{2.30}
where $b$ is an arbitrary complex constant. In this case we have
\beq \BT\BH=b\BV^{-1}\BGD(\Bk)\Bw-b\Bw=-b\BGG(\Bk)\BV\Bw=\widehat{\BL}^\dag(\Bk)\BG,
\eeq{2.31}
where
\beq \BG=-b(\widehat{\BL}^\dag(\Bk)\BV\widehat{\BL}(\Bk))^{-1}\widehat{\BL}^\dag(\Bk)\BV\Bw.
\eeq{2.32}
So the potential $\dund{\BU}^1(\Bx)$ and the $Q^*$-special fields  $\dund{\BE}(\Bx)$ and $\dund{\BJ}(\Bx)$ take the form
\beqa \dund{\BU}^1(\Bx) & =  & -\sum_{\Bk\ne 0}b(\Bk)[(\widehat{\BL}^\dag(\Bk)\BV\widehat{\BL}(\Bk))^{-1}\widehat{\BL}^\dag(\Bk)]\BV\Bw\,e^{i\Bk\cdot\Bx}, \nonum
\dund{\BE}(\Bx) & =  & \BT\dund{\BJ}(\Bx)=\BT\dund{\BJ}^0-\sum_{\Bk\ne 0}b(\Bk)\BGG(\Bk)\BV\Bw\,e^{i\Bk\cdot\Bx}, \nonum
\dund{\BJ}(\Bx) & =  & \dund{\BJ}^0+\sum_{\Bk\ne 0}b(\Bk)\BGD(\Bk)\Bw\,e^{i\Bk\cdot\Bx},
\eeqa{2.33}
where $\BT\dund{\BJ}^0\in \CE^0$ and 
\beq b(\Bk)=0\quad{\rm if}\,\,\Bk\notin\CM.
\eeq{2.33a}

We next look for the polynomial $q(\Bk)$ of lowest degree such that 
\beq \BP(\Bk)=-q(\Bk)(\widehat{\BL}^\dag(\Bk)\BV\widehat{\BL}(\Bk))^{-1}\widehat{\BL}^\dag(\Bk)\BV\Bw
\eeq{2.34}
is also polynomial in $\Bk$. Set
\beq b(\Bk)=q(\Bk)\widehat{c}(\Bk),
\eeq{2.36}
where the $\widehat{c}(\Bk)$ are the Fourier components of some scalar periodic function $c(\Bx)$ such that  now $\widehat{c}(\Bk)$ is zero if $\Bk\notin\CM$.
Then we see that
\beq \dund{\BU}^1(\Bx)=\BP(\Md/\Md x_1, \Md/\Md x_2,\ldots , \Md/\Md x_d)c(\Bx),
\eeq{2.37}
which also gives the $Q^*$-special fields $\dund{\BE}=\BT\dund{\BJ}^0+\BL\dund{\BU}^1$ and $\dund{\BJ}=\BT^{-1}\dund{\BE}$ (in the case
where $\BT$ is non-singular).

The derivation still goes through if the quadratic form associated with $\BV^{-1}$ is $Q^*$-convex or marginally $Q^*$-convex. In the latter case 
one needs to choose $\Bw$ so that 
\beq \BGG_2(\Bk)\Bw\quad{\rm lies\,\,in\,\,the\,\,range\,\,of}\,\,\BGG_2(\Bk)\BV^{-1}\BGG_2(\Bk),\quad \forall\Bk\ne 0,
\eeq{2.38}
since otherwise, from \eq{2.25ab}, the supremum in \eq{2.29} is surely infinite. Thus by induction, successively
setting
\beq \BV^{-1}=\BW_0-\sum_{i=1}^{j-1}\Bw_i\Bw_i^\dag/\Gb_i, \quad\Bw=\Bw_{j},\quad \Gb=\Gb_{j},
\eeq{2.39}
for $j=1,2,\ldots, n$ (where $\BW_0$ has an associated $Q^*$-convex quadratic form, $\Bw=\Bw_{j}$ has to be chosen so that \eq{2.25e} is satisfied, and $\Gb=\Gb_{j}$ is given by \eq{2.29}) one finally reaches a point  
 such that the quadratic function $g$ associated with 
\beq \BT=\BW_0-\sum_{i=1}^n\Bw_i\Bw_i^\dag/\Gb_i 
\eeq{2.40}
is an extremal quadratic $Q^*$-convex function in the sense that one cannot subtract any non-zero positive semi-definite tensor from $\BT$ and retain its $Q^*$-convexity.
In this case, with $\BV^{-1}=\BT$ the supremum in \eq{2.29} is infinite for all non-zero $\Bw\in\CC^m$.

\section{A Generalization}
\setcounter{equation}{0}
Following the ideas of Cherkaev and Gibiansky \cite{Cherkaev:1992:ECB} (see also sections 24.1 and 30.4 of \cite{Milton:2002:TOC})
the arguments presented here extend directly to integrals of the form
\beq \int_{\GO}h(\widetilde{\BE}(\Bx),\BJ(\Bx))\,d{\Bx},
\eeq{4.1}
where $h$ is quadratic and the $m$ component field $\BJ(\Bx)$ still satisfies $\BL^\dag\BJ=0$ while the
$n$-component field $\widetilde{\BE}(\Bx)$ derives from a $p$ component potential $\widetilde{\BU}(\Bx)$ according to the relations 
\beq \widetilde{E}_{r}(\Bx)= \sum_{q=1}^{p}\widetilde{L}_{rq}\widetilde{U}_{q}(\Bx),
\eeq{4.2}
for $r=1,2,\ldots{},n$, where $\widetilde{L}_{rq}$ is the differential operator
\beq \widetilde{L}_{rq}=\widetilde{A}_{rq}+\sum_{h=1}^v\sum_{a_1,\ldots{},a_h=1}^{d} \widetilde{A}_{rqh}^{a_1\ldots{}a_h}
{\Md\over\Md x_{a_1}}{\Md\over\Md x_{a_2}}\ldots{\Md\over\Md x_{a_h}}
\eeq{4.3}
of order $v$ in a space of dimension $d$ with real or complex valued constant coefficients $\widetilde{A}_{rq}$ and
$\widetilde{A}_{rqh}^{a_1\ldots{}a_h}$. 

The function $h(\widetilde{\BE},\BJ)$ can be expressed in the form
\beq h(\widetilde{\BE},\BJ)=\pmatrix{\widetilde{\BE}\cr \BJ}\cdot\BM\pmatrix{\widetilde{\BE}\cr \BJ},
\eeq{4.4}
where the matrix $\BM$ is Hermitian. It is chosen to be sharply $Q^*$-convex, so that for periodic fields $\widetilde{\BE}(\Bx)$ and $\BJ(\Bx)$ satisfying the differential constraints one has
\beq \lang h(\widetilde{\BE},\BJ)\rang\geq h(\lang\widetilde{\BE}\rang,\lang\BJ\rang),
\eeq{4.5}
with equality for some non-constant $Q^*$-special fields $\widetilde{\BE}=\widetilde{\underline{\BE}}$ and $\BJ=\dund{\BJ}$.

\section*{Acknowledgements}
The author thanks Hyeonbae Kang for his collaboration on the paper \cite{Kang:2013:BVF3d} which included the example reviewed in section 2, and thanks
him and Marc Briane for helpful comments. In addition he is grateful to Ben Eggleton, CUDOS and the University of Sydney for the provision of 
office space during his visit there and to the
National Science Foundation for support through grant DMS-1211359.

\bibliographystyle{mod-phaip}
\bibliography{tcbook,newref}

\end{document}


\vspace{-1in}
\title{Addendum to ``Sharp inequalities which generalize the divergence theorem--an extension of the notion of quasiconvexity''}

\author{Graeme W. Milton\\
\small{Department of Mathematics, University of Utah, Salt Lake City UT 84112, USA} \\
\small{(milton@math.utah.edu)}}
\date{}
\maketitle
\begin{abstract}
The paper ``Sharp inequalities which generalize the divergence theorem--an extension of the notion of quasiconvexity'' published in Proc. R. Soc. A. 2013, Volume 469, Article 20130048; doi:10.1098/rspa.2013.0048 is clarified.
Notably, much more general boundary conditions are given under which sharp lower bounds on the integrals of certain quadratic functions of the fields can be obtained. As a consequence, strict $Q^*$-convexity is found to be
an appropriate condition to ensure uniqueness of the solutions of a wide class of partial differential equations
with constant coefficients, in a given domain $\GO$ with appropriate boundary conditions.

\end{abstract}
\vskip2mm

\noindent Keywords: Quasiconvexity, $Q^*$-convexity, Null-Lagrangians 
\newline
\vskip 5mm
In the paper \cite{Milton:2013:SIG}, which generalized a result in \cite{Kang:2013:BVF3d},
 sharp lower bounds on the integral 
\beq \int_\GO f(\BE(\Bx))\,d{\Bx}, \eeq{0.5}
where $\Bx$ lies in a $d$-dimensional space, were obtained for sharply $Q^*$-convex quadratic functions 
$f(\BE(\Bx))$ where $\BE$ derives some potential $\BU$, for certain special boundary conditions
on $\BU$ and its derivatives. Specifically, the field
$\BE$ with $m$ components $E_r(\Bx),~r=1,2,\ldots{},m$ derives from a real or complex
potential $\BU(\Bx)$, with $\ell$ components $U_1(\Bx),\ldots{},U_{\ell}(\Bx)$ through the equations
\beq E_{r}(\Bx)= \sum_{q=1}^{\ell}L_{rq}U_{q}(\Bx), \eeq{1.1}
for $r=1,2,\ldots{},m$, where $L_{rq}$ is the differential operator
\beq L_{rq}=A_{rq}+\sum_{h=1}^t\sum_{a_1,\ldots{},a_h=1}^{d} A_{rqh}^{a_1\ldots{}a_h}
{\Md\over\Md x_{a_1}}{\Md\over\Md x_{a_2}}\ldots{\Md\over\Md x_{a_h}}
\eeq{1.2}
of order $t$ in a space of dimension $d$ with real or complex valued constant coefficients $A_{rq}$ and
$A_{rqh}^{a_1\ldots{}a_h}$ , some of which may be zero. Here the function  $f(\BE(\Bx))$ is said to be 
a $Q^*$-convex function if it satisfies 
\beq \lang f(\BE) \rang \geq f(\lang\BE\rang),
\eeq{1.2a}
for all periodic functions $\BE=\BL\BU$ that derive from a potential $\BU(\Bx)$ that is the sum of a polynomial $\BU^0(\Bx)$ and a periodic potential $\BU^1(\Bx)$ .
The angular brackets denote volume averages over the unit cell of periodicity. (By all periodic functions we also meant for all primitive unit cells of periodicity, including parallelepiped shaped ones).  
The function $f$ is sharply  $Q^*$-convex  if in addition one has the equality
\beq \lang f(\underline{\BE}) \rang = f(\lang\underline{\BE}\rang),
\eeq{1.2b}
for some non-constant periodic function $\underline{\BE}=\BL\underline{\BU}$ that derives from a potential $\underline\BU(\Bx)$ that is the sum of a polynomial $\underline{\BU}^0(\Bx)$ and a periodic potential $\underline{\BU}^1(\Bx)$ . The fields $\underline{\BE}$ and its potential $\underline\BU(\Bx)$ 
were called $Q^*$-special fields.  When the elements $A_{rq}$ vanish and the coefficients  $A_{rqh}^{a_1\ldots{}a_h}$ are zero for all $h\ne t$
then $Q^*$-convexity may be equivalent to quasiconvexity, but otherwise they are not equivalent: see \cite{Ball:1981:NLW} for a precise definition of quasiconvexity.

In fact such inequalities can easily be obtained for many more boundary 
conditions, even for functions that are $Q^*$-convex but not sharply  $Q^*$-convex. An obvious example is when
$\BE=\Grad U$ where $U$ is a scalar, and when $f(\BE)=\BE\cdot\BE$. If we find some potential 
$\underline{U}$ such that $\Delta \underline{U}=0$, then the standard ``energy minimization''
variational principle implies we have the sharp inequality
\beq \int_\GO f(\Grad U)\,d{\Bx}\geq \int_{\Md\GO} \underline{U}(\Grad\underline{U}\cdot\Bn)\,d{S},
\eeq{1.2aa}
for all potentials $U(\Bx)$ satisfying the boundary condition that $U(\Bx)=\underline{U}(\Bx)$ when $\Bx\in\Md\Omega$, where $\Bn$ is the
outward normal to the surface $\Md\GO$ of $\GO$. The inequality is of course sharp when 
$U(\Bx)=\underline{U}(\Bx)$ for all $\Bx$ in $\Omega$. 

Generalizing this, we should define a solution field in $\GO$ to be any field $\underline{\BE}$ satisfying
\beq \underline{\BJ}(\Bx)=\BS\underline{\BE}(\Bx),\quad \underline{\BE}=\BL\underline{\BU},\quad \BL^\dag\underline{\BJ}=0,
\eeq{2.7}
where $\BS$ is the Hermitian $m\times m$ matrix defining the quadratic form $f$: 
\beq f(\BE(\Bx)) =  \BE\cdot\BS\BE\equiv \sum_{r=1}^m\sum_{s=1}^m \overline{E}_r(\Bx)S_{rs}E_s(\Bx), \nonum
\eeq{1.5}
and where $\BL^\dag$ is the formal adjoint of $\BL$, 
\beq L_{qr}^\dag=\overline{A_{rq}}+\sum_{h=1}^t\sum_{a_1,\ldots{},a_h=1}^{d} (-1)^h \overline{A_{rqh}^{a_1\ldots{}a_h}}
{\Md\over\Md x_{a_1}}{\Md\over\Md x_{a_2}}\ldots{\Md\over\Md x_{a_h}},
\eeq{1.4a}
in which the bar denotes complex conjugation. Since the operators $\BL$ and $\BL^\dag$ are formal adjoints the quantity
\beq f_0=\int_\GO f(\underline{\BE}(\Bx))\,d{\Bx}=\int_\GO \underline{\BE}(\Bx)\cdot\BS\underline{\BE}(\Bx)\,d\Bx=\int_\GO \underline{\BE}(\Bx)\cdot\underline{\BJ}(\Bx)\,d\Bx \eeq{3.1}
can be computed in terms of boundary terms using integration by parts. Now we show that 
\beq \int_\GO f(\BE(\Bx))\,d{\Bx}\geq f_0,
\eeq{3.2}
for all fields $\BE(\Bx)$ deriving from a potential $\BU(\Bx)$ that matches the appropriate boundary data of the potential $\underline{\BU}(\Bx)$ (generally involving
both $\BU(\Bx)$ and its derivatives when $t>1$), provided certain further supplementary conditions hold.

To establish the inequality in \eq{3.2} we first find a parallelepiped $C$ that contains $\GO$. Inside $\GO$ define $\Gd\BE(\Bx)=\BE(\Bx)-\underline{\BE}(\Bx)$
and $\Gd\BU(\Bx)=\BU(\Bx)-\underline{\BU}(\Bx)$.
We extend $\Gd\BE(\Bx)$ and $\Gd\BU(\Bx)$ outside $\GO$ so that it is periodic with $C$ as a unit cell. In this cell, but outside $\GO$, we set
$\Gd\BE(\Bx)=0$ and $\Gd\BU(\Bx)=0$. The boundary data on $\GO$ for 
the potential $\BU(\Bx)$ are chosen so the equation $\Gd\BE=\BL\Gd\BU$ holds weakly across the boundary of $\GO$.
Defined in this way, the relation $\Gd\BE=\BL\Gd\BU$ holds in a weak sense, and
so the inequality \eq{1.2a} implies
\beq  \int_\GO f(\Gd\BE)\,d\Bx\geq |C| f(\lang\Gd\BE\rang),
\eeq{3.4}
where $|C|$ is the volume of $C$. Since $f$ is quadratic,
\beq f(\Gd\BE)=f(\BE)-f(\underline\BE)-2\Gd\BE\cdot\BS\underline\BE, \eeq{3.9aa}
and integration by parts implies
\beq
\int_\GO \Gd\BE(\Bx)\cdot\BS\underline{\BE}(\Bx)\,d\Bx=\int_\GO \Gd\BE(\Bx)\cdot\underline{\BJ}(\Bx)\,d\Bx=0.
\eeq{3.9ab}
So it follows from this and \eq{3.4} that
\beq \int_\GO f(\BE)\,d\Bx\geq f_0+ |C| f(\lang\Gd\BE\rang).
\eeq{3.9ac}
This inequality in general requires us to know $\lang\Gd\BE\rang$. Given a $m$-dimensional constant vector $\BJ^0$, we have
\beq \int_C\BJ^0\cdot\Gd\BE\,d\Bx=\int_C\BJ^0\cdot\BL\Gd\BU\,d\Bx=
\int_C\BJ^0\cdot\BA\Gd\BU\,d\Bx,
\eeq{3.7}
where $\BA$ is the matrix with elements $A_{rq}$. (The last equality in \eq{3.7} is established using integration by parts and the fact that
$\Gd\BU$ is zero in the vicinity of the boundary of $C$). Since this holds for all $\BJ^0$, we deduce that
\beq \lang\Gd\BE\rang=\frac{1}{|C|}\int_C \BA\Gd\BU\,d\Bx.
\eeq{3.8}
Therefore a sufficient condition for $f(\lang\Gd\BE\rang)$ to be zero is that
\beq \BS\BA=0, \eeq{3.9}
so that the range of $\BA$ is in the null space of $\BS$ (which if $\BS$ is non-singular requires that $\BA=0$). When this supplementary condition holds then clearly
\eq{3.9ac} implies the inequality in \eq{3.2}.

One corollary of this analysis is that strict $Q^*$-convexity is an appropriate condition to ensure uniqueness. 
Strict $Q^*$-convexity means that if \eq{3.2} is satisfied as an equality, then $\BE(\Bx)=\lang\BE\rang$.
If there were two solution fields $\underline\BE$ and $\underline\BE'$ sharing the same boundary data  on $\GO$, 
and the supplementary condition \eq{3.9} held, then the inequality \eq{3.2} would imply
\beq  \int_\GO f(\underline\BE')\,d\Bx\geq \int_\GO f(\underline\BE)\,d\Bx. \eeq{3.10a}
By interchanging $\underline\BE$ and $\underline\BE'$, the reverse inequality must also hold. So one
must have equality. Tracing back the argument which lead to \eq{3.2} one must also have equality in
\eq{3.4} when inside $\GO$, $\BE$ is replaced by $\underline\BE'$, which by the strict $Q^*$-convexity implies
\beq \Gd\BE(\Bx)=\lang\Gd\BE\rang. \eeq{3.10b}
As the left hand side is zero in $C\setminus\GO$, the right hand side must be zero also and we 
conclude that inside $\GO$,
\beq \Gd\BE(\Bx)=\underline\BE'(\Bx)-\underline\BE(\Bx)=0. \eeq{3.10c}
Even if $f$ is sharply $Q^*$-convex, rather than strictly $Q^*$-convex, uniqueness could still hold.
Indeed the field $\Gd\BE(\Bx)$ must be a $Q^*$-special field which in addition vanishes in $C\setminus\GO$,
and there may be no non-zero such fields.

Examples of the solution fields $\underline\BE(\Bx)$ are easily found. Following, for example section 6.9
in \cite{Milton:2002:TOC}, we could look
for fields where the potentials take the form
\beq \underline\BU(\Bx)=\Bu h(i\Bx\cdot\Bk), \eeq{3.10d}
where $\Bu$ is a complex $m$-component vector, $\Bk$ is a complex $d$ dimensional vector, and $h(z)$ is 
an analytic function of $z$ such that the singularities of $h(\Bx\cdot\Bc)$ lie outside $\GO$. The 
function $h(z)$ could be chosen to be $e^z$ in which case $\underline\BE(\Bx)$ is a plane wave but it is
not necessary to do so: we could, for example, take $h(z)=1/(z-z_0)$. Introduce the elements
\beqa 
\widehat{L}_{rq}(\Bk) & = & A_{rq}+\sum_{h=1}^t\sum_{a_1,\ldots{},a_h=1}^{d} i^hA_{rqh}^{a_1\ldots{}a_h}
k_{a_1}k_{a_2}\ldots{}k_{a_h}, \nonum
\widehat{L}_{qr}^\dag (\Bk) & = & \overline{A_{rq}}+\sum_{h=1}^t\sum_{a_1,\ldots{},a_h=1}^{d} (-i)^h\overline{A_{rqh}^{a_1\ldots{}a_h}}
k_{a_1}k_{a_2}\ldots{}k_{a_h},
\eeqa{1.7}
of the $m\times \ell$ matrix $\widehat{\BL}(\Bk)$ and its adjoint $\widehat{\BL}^\dag (\Bk)$, 
defined for all $\Bk\in\CC^d$. Then $\underline\BE(\Bx)$ will be a solution field if and only if
\beq \BM(\Bk)\Be=0,\quad{\rm where}~\BM(\Bk)=\widehat{\BL}^\dag (\Bk)\BS\widehat{\BL}(\Bk).
\eeq{3.10e}
In other words, we should choose the complex-valued vector $\Bk$ to be a solution of
\beq \det[\BM(\Bk)]=0, \eeq{3.10f}
and then take $\Be$ to be an associated null-vector of $\BM(\Bk)$. More generally, a wide variety
of solution fields could be generated by superpositions of fields of the form \eq{3.10d}. 

Similarly we can bound integrals of the form 
\beq  \int_\GO g(\BJ(\Bx))\,d{\Bx}, \eeq{1.3}
for fields $\BJ(\Bx)$ with $m$ real or complex components $J_r(\Bx),~r=1,2,\ldots{},m$,  that satisfy the differential constraints 
\beq \sum_{r=1}^{m}L_{qr}^\dag J_{r}(\Bx)=0,
\eeq{1.4} 
and for quadratic functions $g(\BJ)$ which are $Q^*$-convex in the sense that the inequality
\beq \lang g(\BJ) \rang \geq g(\lang\BJ\rang),
\eeq{1.4b}
holds for all periodic functions $\BJ$ satisfying $\BL^\dag\BJ=0$. 

In this context we should define the solution
fields $\dund{\BJ}$ in $\GO$ to be those fields for which there exists at least one potential $\dund{\BU}$ 
in $\GO$ such that 
\beq \dund{\BE}(\Bx)=\BT\dund{\BJ}(\Bx),\quad \dund{\BE}=\BL\dund{\BU},\quad \BL^\dag\dund{\BJ}=0,
\eeq{2.14}
where $\BT$ is the Hermitian $m\times m$ matrix entering the quadratic form $g$:
\beq g(\BJ(\Bx)) =  \BJ\cdot\BT\BJ\equiv \sum_{r=1}^m\sum_{s=1}^m \overline{J}_r(\Bx)T_{rs}J_s(\Bx).
\eeq{1.5a}
The quantity 
\beq  g_0=\int_\GO g(\dund{\BE}(\BJ))\,d{\Bx}=\int_\GO \dund{\BJ}(\Bx)\cdot\BT\dund{\BJ}(\Bx)=\int_\GO \dund{\BJ}(\Bx)\cdot\dund{\BE}(\Bx)
\eeq{3.1b}
can be computed in terms of boundary values using integration by parts. 
Define $\Gd\BJ(\Bx)=\BJ(\Bx)-\dund{\BJ}(\Bx)$. We extend $\Gd\BJ(\Bx)$ outside $\GO$ so that it is are periodic with $C$ as a unit cell. In this cell, but outside $\GO$, 
$\Gd\BJ(\Bx)=0$. The boundary data on $\GO$ for $\BJ(\Bx)$ are chosen so the equation $\BL^\dag\Gd\BJ=0$
holds weakly across the boundary of $\GO$. Defined in this way, the relation $\BL^\dag\Gd\BJ=0$ 
holds in a weak sense, and analogously to \eq{3.9ac} we have
\beq  \int_\GO g(\BJ)\,d\Bx\geq g_0+ |C|g(\lang\Gd\BJ\rang).
\eeq{3.12}
This inequality in general requires us to know $\lang\Gd\BJ\rang$. Let $\CE^0$ 
denote the vector space of $m$ dimensional constant vectors $\BE^0$ expressible in the form
$\BE^0=\BL\BU^0(\Bx)$ for some polynomial potential $\BU^0(\Bx)$. Given such a vector $\BE^0$, we have
\beq \int_C \Gd\BJ\cdot{\BE^0}\,d\Bx=\int_C \Gd\BJ\cdot\BL{\BU^0}\,d\Bx=0,
\eeq{3.13}
where the last equality follows from integration by parts, using the fact that $\Gd\BJ=0$ in a vicinity of the boundary of $C$. Since this holds for all $\BE^0\in\CE^0$ we deduce that 
a sufficient condition for $g(\lang\Gd\BJ\rang)$ to be equal to be zero is that
the range of $\BT$ be a subset of $\CE^0$ and if this ancillary condition holds then clearly \eq{3.12} implies 
\beq  \int_\GO g(\BJ)\,d\Bx\geq g_0.  \eeq{3.13a}

Again strict $Q^*$-convexity is an appropriate condition to ensure uniqueness. 
Strict $Q^*$-convexity means that if \eq{1.4b} is satisfied as an equality, then $\BJ(\Bx)=\lang\BJ\rang$.
If there were two solution fields $\dund\BJ$ and $\dund\BJ'$ sharing the same boundary data on $\Md\GO$, 
and the ancillary condition that the range of $\BT$ be a subset of $\CE^0$ is satisfied,
then \eq{3.13a} holds as an equality when $\BJ=\dund\BJ'$. This in turn implies $\dund\BJ=\dund\BJ'$
inside $\GO$. 

In closing we remark that for all the results to hold it clearly suffices for $f$ or $g$ to be 
$Q^*_C$-convex, where by $Q^*_C$-convex we mean that the inequalities \eq{1.2a} or
\eq{1.4b} hold for $C$-periodic functions with a unit cell $C$ containing $\GO$, and not
necessarily for all periodic functions.

\section*{Acknowledgements}
The author thanks the National Science Foundation for support through grant DMS-1211359.

\bibliography{/home/milton/tcbook,/home/milton/newref}
\bibliographystyle{/home/milton/latex/mod-phaip}